# NONPARAMETRIC ESTIMATION WHEN DATA ON DERIVATIVES ARE AVAILABLE

By Peter Hall and Adonis Yatchew[1]

*Australian National University and University of Toronto*

We consider settings where data are available on a nonparametric function and various partial derivatives. Such circumstances arise in practice, for example in the joint estimation of cost and input functions in economics. We show that when derivative data are available, local averages can be replaced in certain dimensions by nonlocal averages, thus reducing the nonparametric dimension of the problem. We derive optimal rates of convergence and conditions under which dimension reduction is achieved. Kernel estimators and their properties are analyzed, although other estimators, such as local polynomial, spline and nonparametric least squares, may also be used. Simulations and an application to the estimation of electricity distribution costs are included.

**1. Introduction.** We consider settings where data are available on a nonparametric function and various partial derivatives. For example, suppose data $(X_{1i}, X_{2i}, Y_i, Y_{1i})$, $i = 1, \ldots, n$, are available for

$$y = g(x_1, x_2) + \varepsilon, \qquad y_1 = \frac{\partial g(x_1, x_2)}{\partial x_1} + \varepsilon_1.$$

Then $g$ can be estimated at rates as though it were a function of a single nonparametric variable, rather than two. Heuristically, the presence of data on the partial derivative with respect to $x_1$ eliminates the need for local averaging in the $x_1$ direction. This, in turn, results in dimension reduction and suggests the possibility of estimating $g$ and its derivatives at relatively fast rates.

Received March 2005; revised November 2005.
[1]Supported in part by a grant from the Social Sciences and Humanities Research Council of Canada.
*AMS 2000 subject classifications.* Primary 62G08, 62G20; secondary 62P20, 91B38.
*Key words and phrases.* Dimension reduction, kernel methods, nonparametric regression, partial derivative data, rates of convergence, statistical smoothing, cost function estimation.

It is natural to ask whether data on derivatives would be available in practical settings, or whether this investigation is esoteric. In fact, such data are commonly available in economics. The underlying reason is that economic models frequently assume that agents economize, that is, that they implicitly or explicitly solve constrained optimization problems. Thus, data may not only be available on an objective function, but also on first order conditions related to the optimization problem. An example serves to illustrate the point.

Consider $y = g(Q, r, w) + \varepsilon$, where $y$ is the minimum cost of producing output level $Q$ given $r$ and $w$, the prices of capital and labor, respectively, and $\varepsilon$ is a residual. By the envelope theorem or, equivalently, Shephard's Lemma (see, e.g., [17, 23]), the partial derivatives of $g$ with respect to $r$ and $w$ yield the optimal quantities of capital and labor required to produce $Q$. Joint estimation of cost functions and their partial derivatives (i.e., the inputs) using parametric models is routinely undertaken (see, e.g., [10, 14]). Florens, Ivaldi and Larribeau [7] analyze the behavior of parametric approximations of systems such as the ones considered in this paper. However, nonparametric estimation as proposed here has received little attention.

Quite different examples of the same type arise in engineering settings, for example, in real-time records of certain types of motion sensors and in modeling problems connected to stochastic control; see, for instance, [16, 19].

Rates of convergence for nonparametric regression (e.g., [20, 21]) often limit the usefulness of conventional nonparametric models in fields where regression modeling involves multiple explanatory variables. Several devices are available to mitigate this curse of dimensionality. They include additive and varying-coefficient models (see, e.g., [3, 11, 12, 22]), projection-based methods (e.g., [4, 9, 13]), and recursive partitioning and tree-based methods (e.g., [2, 8, 26]). For the most part, these approaches fit "abbreviated" models, where components or interactions among components are dropped in order to reduce the variability of an estimator. We shall show that incorporating derivative information can yield lower variability and faster convergence rates for the full underlying regression function, without any need for abbreviation.

Methodology based on this idea can be expected to reach beyond examples in economics and engineering such as those given earlier. Particularly with the development of new technologies which allow rates of change to be recorded at discrete times, systems in the physical and biological sciences offer opportunities for dimension reduction using derivative data. For example, in meteorology, each of barometric pressure, wind speed and direction (the latter two being functions of the pressure gradient) are measured over broad geographic regions. In some fields, an evolving system is often modeled as a (possibly constrained) optimization problem, so one might expect data relating to first order conditions to be available there.

This paper is organized as follows. Section 2 outlines our assumptions and provides results on optimal rates of convergence. The approaches to dimension reduction addressed there are nonstandard. Section 3, which shows that suitably constructed kernel estimators achieve optimal rates of convergence, uses familiar smoothing methods surveyed by, for example, Wand and Jones [25], Fan and Gijbels [5] and Simonoff [18]. We also note that the idea of combining local and nonlocal averaging has been used by Linton and Nielsen [15] and Fan, Härdle and Mammen [6]. Results of Bickel and Ritov [1] on estimators that are constructed by "plugging in" root-$n$ consistent estimators of functions are more distantly related. Section 4 describes results of simulations and an empirical application involving data on electricity distribution costs. Proofs of propositions are deferred to the Appendix.

Before proceeding, it may be useful to illustrate our results on rates of convergence. Let $g(x_1, x_2, x_3, x_4)$ be a nonparametric function for which we have data and consider the following hierarchy of functions where superscripts denote multiple first order partial derivatives:

$$
\begin{array}{ccccccc}
& & & g^{(1,1,1,1)} & & & \\
& & g^{(1,1,1,0)} & g^{(1,1,0,1)} & g^{(1,0,1,1)} & g^{(0,1,1,1)} \\
g^{(1,1,0,0)} & g^{(1,0,1,0)} & g^{(0,1,1,0)} & g^{(1,0,0,1)} & g^{(0,1,0,1)} & g^{(0,0,1,1)} \\
g^{(1,0,0,0)} & g^{(0,1,0,0)} & g^{(0,0,1,0)} & g^{(0,0,0,1)} & & \\
\end{array}
$$

(a) If data are available on the complete hierarchy, then $g$ can be estimated root-$n$ consistently—that is, the "nonparametric dimension" of the estimation problem is zero. (b) If data are available on all multiple first order partials for any subset of $p$ variables, then the nonparametric dimension is $4 - p$. For example, if one observes all partials below the main diagonal, then the nonparametric dimension is one. (c) If data are available on all multiple first order partials for any subset of $p$ variables, except those appearing in the bottom $\bar{\ell}$ rows, then the nonparametric dimension is $4 - (p - \bar{\ell})$. For example, if one observes all partials in the northwest wedge, then the nonparametric dimension is two. (d) For an arbitrary set of observed partial derivatives, an upper bound on the nonparametric dimension of the estimation problem may be determined by using (b) and (c) to find the subset which yields the lowest nonparametric dimension. For example, if one observes all simple first order partials, that is, all partials in the bottom row, then the nonparametric dimension does not exceed three. If, in addition, one observes $g^{(1,1,0,0)}$, then the nonparametric dimension does not exceed two.

## 2. Properties underpinning methodology.

2.1. *Main theorem about functionals.* For simplicity, we shall assume that $g$ is supported on the unit cube $\mathcal{R}_k = [0,1]^k$, although substantially more general designs are possible. Let $A$ denote the set of all sequences $\alpha = (\alpha_1, \ldots, \alpha_k)$ of length $k$ consisting solely of zeros and ones. Given $\alpha \in A$ and $x = (x_1, \ldots, x_k) \in \mathcal{R}_k$, define $|\alpha| = \sum_j \alpha_j$ and

$$g^\alpha(x) = \frac{\partial^{|\alpha|} g(x)}{\partial x_{i_1} \cdots \partial x_{i_{|\alpha|}}},$$

where $i_1 < \cdots < i_{|\alpha|}$ denotes the sequence of indices $i$ for which $\alpha_i = 1$.

Let $\mathcal{B}_k$ denote the class of bounded functions on $\mathcal{R}_k$ and let $\mathcal{G}_k$ denote the class of functions $g$ on $\mathcal{R}_k$ for which $g^\alpha \in \mathcal{B}_k$ for each $\alpha \in A$. Given $C > 0$, let $\mathcal{K}(C)$ denote the class of functionals $\psi$ that may be represented as

$$(2.1) \qquad (\psi g)(x) = \int_{\mathcal{R}_k} \chi(u, x) \, g(u) \, du \qquad \text{for all } g \in \mathcal{B}_k,$$

where the function $\chi$ (which determines $\psi$) satisfies $\sup_{u,x \in \mathcal{R}_k} |\chi(u,x)| \leq C$.

THEOREM 1. *There exists a set of functionals $\{\psi_\alpha, \alpha \in A\} \subseteq \mathcal{K}(1)$ such that for all $g \in \mathcal{G}_k$,*

$$g = \sum_{\alpha \in A} \psi_\alpha g^\alpha.$$

A proof of this theorem and explicit formulae for the functionals $\psi_\alpha$ are given in Appendix A.1.

To appreciate the implications of Theorem 1 for inference, assume that for each $\alpha \in A$, that is, for each model $y_\alpha = g^\alpha(x) + \varepsilon_\alpha$, we have data pairs $(X_{\alpha i}, Y_{\alpha i})$ generated by

$$(2.2) \qquad Y_{\alpha i} = g^\alpha(X_{\alpha i}) + \varepsilon_{\alpha i}, \qquad 1 \leq i \leq n_\alpha,$$

where the $X_{\alpha i}$'s are distributed on $\mathcal{R}_k$ with a density $f_\alpha$ that is bounded away from zero there and the errors $\varepsilon_{\alpha i}$ are independent with zero means and bounded variances, also independent of the $X_{\alpha i}$'s. Put $n = \min_{\alpha \in A} n_\alpha$. It follows from the form of the functional $\psi_\alpha$ [see (2.1)] that from these data, we may construct an estimator $\widehat{\psi_\alpha g^\alpha}$ of $\psi_\alpha g^\alpha$ that is root-$n$ consistent whenever $g \in \mathcal{G}_k$.

For example, if the $X_{\alpha i}$'s are uniformly distributed on $\mathcal{R}_k$ and if $\psi = \psi_\alpha$ is given by (2.1) with $\chi$ there denoted by $\chi_\alpha$, then an unbiased root-$n$ consistent estimator of $\psi_\alpha g^\alpha$ is given by $\widehat{\psi_\alpha g^\alpha}$, where

$$(2.3) \qquad (\widehat{\psi_\alpha g^\alpha})(x) = \frac{1}{n_\alpha} \sum_{i=1}^{n_\alpha} Y_{\alpha i} \chi_\alpha(X_{\alpha i}, x).$$

Theorem 1 now implies that

$$\hat{g} \equiv \sum_{\alpha \in A} \widehat{\psi_\alpha g^\alpha} \tag{2.4}$$

is a root-$n$ consistent estimator of $g$. Properties of estimators such as $\widehat{\psi_\alpha g^\alpha}$ and $\hat{g}$ will be discussed in Section 3.2.

The theory that we develop admittedly does not address the "cost" of sampling data on derivatives. In the examples from economics and engineering discussed in Section 1, the cost is low, although in some other problems it is prohibitively high. Moreover, if high order derivative information is absent, then our estimators simply do not enjoy fast convergence rates. We characterize convergence rates in terms of the value of $n = \min_{\beta \in B} n_\beta$ and do not dwell on the fact that if there is a sufficiently large order of magnitude of data on $(X, Y)$ alone, sufficiently greater than $n$, then the convergence rate of a conventional nonparametric estimator based solely on those data can be faster than the rates given in this paper.

The assumption that errors for measurements of different derivatives are independent can be significantly relaxed without affecting the theoretical results that we shall give in Section 3. The assumption may not be completely plausible in the setting of capital and labor costs, but it is realistic in the context of engineering problems, where motion sensor data on functions and their derivatives are estimated by different sensors with different characteristics. Correlations among the errors for different functions will be permitted in the simulation study in Section 4.

In the following examples, the decomposition of $g$ provided in Theorem 1 is rearranged to illustrate root-$n$ consistent estimation.

EXAMPLE 1. Suppose $k = 1$ and that noisy data are available for $g(x)$ and $dg(x)/dx$. Write $g(x) \equiv \underline{g}_1(x) + \underline{g}_0(\cdot)$, where

$$\underline{g}_1(x) \equiv \int_0^x \frac{dg(u)}{du} du = g(x) - g(0),$$

$$\underline{g}_0(\cdot) \equiv \int_0^1 \{g(u) - \underline{g}_1(u)\} du = g(0).$$

The function $\underline{g}_1$ can be estimated root-$n$ consistently, in which case its integral and hence $\underline{g}_0$ can too.

EXAMPLE 2. Suppose $k = 2$ and that noisy data are available for $g^{(1,1)}$, $g^{(1,0)}$, $g^{(0,1)}$ and $g^{(0,0)} = g$. Write $g(x) \equiv \underline{g}_{11}(x_1, x_2) + \underline{g}_{10}(x_1, \cdot) + \underline{g}_{01}(\cdot, x_2) + \underline{g}_{00}(\cdot, \cdot)$, where

$$\underline{g}_{11}(x_1, x_2) \equiv \int_0^{x_1} \int_0^{x_2} g^{(1,1)}(u_1, u_2) \, du_1 \, du_2$$

$$= g(x_1, x_2) - g(x_1, 0) - g(0, x_2) + g(0, 0),$$

$$\underline{g}_{10}(x_1, \cdot) \equiv \int_0^{x_1} \int_0^1 g^{(1,0)}(u_1, x_2) \, du_1 \, dx_2 - \int_0^1 \underline{g}_{11}(x_1, x_2) \, dx_2$$

$$= g(x_1, 0) - g(0, 0),$$

$$\underline{g}_{01}(\cdot, x_2) \equiv \int_0^1 \int_0^{x_2} g^{(0,1)}(x_1, u_2) \, dx_1 \, du_2 - \int_0^1 \underline{g}_{11}(x_1, x_2) \, dx_1$$

$$= g(0, x_2) - g(0, 0),$$

$$\underline{g}_{00}(\cdot, \cdot) \equiv \int_0^1 \int_0^1 g(x_1, x_2) - \underline{g}_{11}(x_1, x_2) - \underline{g}_{10}(x_1, \cdot) - \underline{g}_{01}(\cdot, x_2) \, dx_1 \, dx_2$$

$$= g(0, 0).$$

Sample analogues of all integral expressions can be calculated without local averaging. Thus, $\underline{g}_{11}$ and its integrals can be estimated root-$n$ consistently, in which case $\underline{g}_{10}$ and $\underline{g}_{01}$, their respective integrals and $\underline{g}_{00}$ can too.

2.2. *Application of Theorem 1 to lower-dimensional structures.* Let $1 \leq p \leq k$ and consider a lower-dimensional "subspace" of $A$, specifically the set $B$ of all sequences $\beta = (\beta_1, \ldots, \beta_p)$ of length $p$ consisting solely of zeros and ones. Given $g \in \mathcal{B}_k$, define $|\beta|$ and $g^\beta$ analogously to $|\alpha|$ and $g^\alpha$. In particular, $g^\beta$ is a function on $\mathcal{R}_k$, not on the lower-dimensional space $\mathcal{R}_p = [0,1]^p$ and

$$(2.5) \qquad g^\beta(x) = \frac{\partial^{|\beta|} g(x)}{\partial x_{i_1} \cdots \partial x_{i_{|\beta|}}},$$

where $i_1 < \cdots < i_{|\beta|}$ denotes the sequence of indices $i$ for which $\beta_i = 1$. Similarly, although the functional $\psi_\beta$ (the $p$-dimensional analogue of $\psi_\alpha$ introduced in Theorem 1) would normally be interpreted as the functional which takes $b \in \mathcal{B}_p$ to $\psi_\beta b$, defined by

$$(\psi_\beta b)(x_1, \ldots, x_p) = \int_{\mathcal{R}_p} \chi(u_1, \ldots, u_p, x_1, \ldots, x_p) \, b(u_1, \ldots, u_p) \, du_1 \cdots du_p,$$

it can just as easily be interpreted as the functional that takes $g \in \mathcal{B}_k$ to $\psi_\beta g$, defined by

$$(\psi_\beta g)(x_1, \ldots, x_k) = \int_{\mathcal{R}_p} \chi(u_1, \ldots, u_p, x_1, \ldots, x_p)$$

$$(2.6) \qquad \qquad \times g(u_1, \ldots, u_p, x_{p+1}, \ldots, x_k) \, du_1 \cdots du_p.$$

We shall adopt the latter interpretation.

We may, of course, interpret $\beta$ as a $k$-vector and an element of $A$, with its last $k - p$ components equal to zero. We shall take this view in Section 2.3, where we shall treat cases that cannot be readily subsumed under a model in which noisy observations are made of $\psi_\beta g^\beta$ for each $\beta \in B$.

Let $\mathcal{G}_{kp}$ denote the class of functions $g \in \mathcal{B}_k$ for which $g^\beta$ is well defined and bounded on $\mathcal{R}_k$ for each $\beta \in B$. The following result is an immediate corollary of Theorem 1. It is derived by applying Theorem 1 to the function that is defined on $\mathcal{R}_p$ and is obtained from $g$ by fixing the last $k-p$ coordinates of $x$ and allowing the first $p$ coordinates to vary in $\mathcal{R}_p$. However, although Corollary 1 can be proved from Theorem 1, the theorem is a special case of the corollary.

COROLLARY 1. *Assume $1 \leq p \leq k$ and let $\psi_\beta$, for $\beta \in B$, denote the functionals introduced in Section 2.1, but interpreted in the sense of (2.6). Then for each $g \in \mathcal{G}_{kp}$,*

$$g = \sum_{\beta \in B} \psi_\beta g^\beta.$$

The main statistical implication of the corollary is that by observing data on $g^\beta$ for each $\beta \in B$, we reduce the effective dimension of the problem of estimating $g$ from $k$ to $k-p$. The manner in which $g$ depends on its first $p$ components can be estimated root-$n$ consistently and then performance in the estimation problem is driven by the difficulty of determining the way in which $g$ is influenced by its last $k-p$ components.

To better appreciate this point, assume that for each $\beta \in B$, data $(X_{\beta i}, Y_{\beta i})$ are generated by an analogue of the model at (2.2),

(2.7) $$Y_{\beta i} = g^\beta(X_{\beta i}) + \varepsilon_{\beta i}, \qquad 1 \leq i \leq n_\beta,$$

where $g \in \mathcal{G}_{kp}$ and the $X_{\beta i}$'s are distributed on $\mathcal{R}_k$. Suppose, for simplicity, that the common density of the $X_{\beta i}$'s is uniform on $\mathcal{R}_k$. Let $X_{\beta i}^{[k-p]}$ and $x^{[k-p]}$ represent the $(k-p)$-vectors comprised of the last $k-p$ components of $X_{\beta i}$ and $x$, respectively. Denote by $K$ a $(k-p)$-dimensional kernel function, let $h$ be a bandwidth and in close analogy with (2.3), define $\widehat{\psi_\beta g^\beta}$ by

(2.8) $$(\widehat{\psi_\beta g^\beta})(x) = \frac{1}{n_\beta h^{k-p}} \sum_{i=1}^{n_\beta} Y_{\beta i} \chi_\beta(X_{\beta i}, x) K\left(\frac{X_{\beta i}^{[k-p]} - x^{[k-p]}}{h}\right).$$

Set $n = \min_{\beta \in B} n_\beta$. It is readily proved that if (i) $g$ has $d$ derivatives of its last $k-p$ components as well as all multiple first derivatives of its first $p$ components, (ii) $K$ is a bounded, compactly supported, $d$th order kernel, (iii) $x$ is an interior point of $\mathcal{R}_k$, so as to avoid edge effects and, (iv) $h = h(n) \sim$ const $\cdot n^{1/(2d+k-p)}$, then $(\widehat{\psi_\beta g^\beta})(x)$ converges to $(\psi_\beta g^\beta)(x)$ at the standard squared-error rate, $n^{-2d/(2d+k-p)}$, for estimating functions of $k-p$ variables with $d$ derivatives. This result is a consequence of the fact that the smoothing at (2.8) is only over the last $k-p$ coordinates of the data $X_{\beta i}$. Therefore, the estimator

$$\hat{g} \equiv \sum_{\beta \in B} \widehat{\psi_\beta g^\beta}, \tag{2.9}$$

analogous to that at (2.4), converges to $g$ at the squared-error rate $n^{-2d/(2d+k-p)}$. Properties of $\widehat{\psi_\beta g^\beta}$ and $\hat{g}$ will be discussed in Section 3.2.

EXAMPLE 3. Returning to the example in the introduction, suppose $k=2$ and that noisy data are available for $g$ and $g^{(1,0)}$. Write $g(x) \equiv \underline{g}_{11}(x_1,x_2) + \underline{g}_{10}(\cdot, x_2)$, where

$$\underline{g}_{11}(x_1, x_2) \equiv \int_0^{x_1} g^{(1,0)}(u_1, x_2)\, du_1 = g(x_1, x_2) - g(0, x_2),$$

$$\underline{g}_{01}(\cdot, x_2) \equiv \int_0^1 \{g(x_1, x_2) - \underline{g}_{11}(x_1, x_2)\}\, dx_1 = g(0, x_2).$$

Estimates of $\underline{g}_{11}$ and $\underline{g}_{01}$ require local averaging in the $x_2$ direction only. Thus, $\underline{g}_{11}$ can be estimated at one-dimensional optimal rates, in which case its integral and $\underline{g}_{01}$ can too.

EXAMPLE 4. Suppose $k=3$ and that noisy data are available for $g^{(1,1,0)}$, $g^{(1,0,0)}$, $g^{(0,1,0)}$ and $g^{(0,0,0)} = g$. Write $g(x) \equiv \underline{g}_{111}(x_1,x_2,x_3) + \underline{g}_{101}(x_1,\cdot,x_3) + \underline{g}_{011}(\cdot,x_2,x_3) + \underline{g}_{001}(\cdot,\cdot,x_3)$, where

$$\underline{g}_{111}(x_1, x_2, x_3) \equiv \int_0^{x_1}\int_0^{x_2} g^{(1,1,0)}(u_1, u_2, x_3)\, du_1\, du_2$$
$$= g(x_1, x_2, x_3) - g(x_1, 0, x_3) - g(0, x_2, x_3) + g(0, 0, x_3),$$

$$\underline{g}_{101}(x_1, \cdot, x_3) \equiv \int_0^{x_1}\int_0^1 g^{(1,0,0)}(u_1, x_2, x_3)\, du_1\, dx_2 - \int_0^1 \underline{g}_{111}(x_1, x_2, x_3)\, dx_2$$
$$= g(x_1, 0, x_3) - g(0, 0, x_3),$$

$$\underline{g}_{011}(\cdot, x_2, x_3) \equiv \int_0^1\int_0^{x_2} g^{(0,1,0)}(x_1, u_2, x_3)\, dx_1\, du_2 - \int_0^1 \underline{g}_{111}(x_1, x_2, x_3)\, dx_1$$
$$= g(0, x_2, x_3) - g(0, 0, x_3),$$

$$\underline{g}_{001}(\cdot, \cdot, x_3) \equiv \int_0^1\int_0^1 \{g(x_1, x_2, x_3) - \underline{g}_{111}(x_1, x_2, x_3) - \underline{g}_{101}(x_1, \cdot, x_3)$$
$$- \underline{g}_{011}(\cdot, x_2, x_3)\}\, dx_1\, dx_2$$
$$= g(0, 0, x_3).$$

Estimates of each of the above component functions require local averaging in the $x_3$ direction only. Thus, $\underline{g}_{111}$ and its integrals can be estimated at one-dimensional optimal rates, as can $\underline{g}_{101}$ and $\underline{g}_{011}$, their respective integrals and hence also $\underline{g}_{001}$.

With a mild abuse of notation, suppose that $x_3$ in Example 4 is of length $k-2$. Then $g$ can be estimated at $(k-2)$-dimensional optimal rates.

2.3. *More general settings.* In Corollary 1, we assumed that we have available all multiple first derivatives $g^\beta$ of the first $p$ components of $g$. Our restriction to the first $p$ components was made only for notational convenience; they could have been any $p$ components. In particular, we may alter the definition at (2.5) to

$$(2.10) \qquad g^\beta(x) = \frac{\partial^{|\beta|} g(x)}{\partial x_{I(i_1)} \cdots \partial x_{I(i_{|\beta|})}},$$

where $I(1) < \cdots < I(p)$ denotes any given subsequence of length $p$ of $1, \ldots, k$, without affecting the validity of the corollary. The functional $\psi_\beta$ would be interpreted analogously. Taking this view (which we shall in the present section), we may interpret $\beta$ as a $k$-vector.

Low-dimensional cases, such as that treated by Corollary 1, are motivated by circumstances where multiple first derivatives are observed for a subset of variables. It may be that one is able to observe data on $g^\alpha$ for all $\alpha \in A$ such that $|\alpha| \geq \ell$, say, but not for any other values of $\alpha$. This case is not immediately covered by Theorem 1 or Corollary 1, which can be viewed as treating the contrary setting $|\alpha| \leq \ell$.

We shall adopt the general setting discussed in the paragraph containing (2.10) so as to stress the wide applicability of our results. Assume $1 \leq p \leq k$, $0 \leq \ell \leq k$ and $1 \leq p - \ell + 1 \leq k$ and suppose that we have derivative information from components in $\mathcal{P} = \{I(1), \ldots, I(p)\}$. Let $\beta$ and $g^\beta$ be as in (2.10) and assume that we have noisy data on $g^\beta$ for all $\beta \in B$ such that $|\beta| \geq \ell$, as well as for $\beta = 0$; see (2.7). Then we may construct an estimator of $g$, closely analogous to that at (2.9) and enjoying the squared-error convergence rate $n^{-2d/(2d+k-q)}$, where $q = p - \ell + 1$. That rate is valid under the assumption that $g$ has $d$ bounded derivatives.

This result is a consequence of Theorem 2 below, for which we now give notation. Given $\alpha \in A$, $u, x \in \mathcal{R}_k$ and a function $b \in \mathcal{B}_k$, let $i_1, \ldots, i_{|\alpha|}$ denote the indices of the components of $\alpha$ that equal 1. Define $v_\alpha(u, x)$ to be the $k$-vector with $u_{i_j}$ in position $i_j$ for $1 \leq j \leq |\alpha|$ and $x_j$ in position $j$ for each $j$ that is not among $i_1, \ldots, i_{|\alpha|}$. Define the operator $M_\alpha$ by

$$(2.11) \qquad (M_\alpha b)(x) = \int_0^1 \cdots \int_0^1 b\{v_\alpha(u, x)\} \, du_{i_1} \cdots du_{i_{|\alpha|}}.$$

Consider the functional that takes $g$ to the function of which the value at $x$ is

$$\int \xi_\alpha(u, x) \, g^\alpha(u) \, du, \tag{2.12}$$

where $\xi_\alpha(u, x)$ is a function of the $2k$ variables among the components of $u$ and $x$. In Appendix A.2, we shall prove the following result.

THEOREM 2. *If $g \in \mathcal{G}_k$, $1 \leq p \leq k$, $0 \leq \ell \leq k$ and $1 \leq p - \ell + 1 \leq k$, then $g$ can be expressed as a linear form in integrals of the type (2.12), where $|\alpha| \geq \ell$, all components of $\alpha$ that equal 1 are indexed in $\mathcal{P}$ and $\sup_{u,x \in \mathcal{R}_k} |\xi_\alpha(u, x)| \leq C$, with $C > 0$ depending only on $k$, $\ell$ and $p$, and in integrals $M_\beta g$, with $\beta \in A$ and $|\beta| \geq p - \ell + 1$.*

Our derivation of Theorem 2 will provide an inductive argument for calculating the representation of $g$ in any given case.

To appreciate how the convergence rate given three paragraphs above follows from Theorem 2, let us consider the case $p = k$, for simplicity, and express $g$ as indicated in the theorem: $g = g_1 + g_2$, where

$$g_1(x) = \sum_{i=1}^{r} \int_{\mathcal{R}_k} \xi_{\alpha(i)}(u, x) \, g^{\alpha(i)}(u) \, du, \qquad g_2(x) = \sum_{i=1}^{s} c_i \, (M_{\beta(i)} g)(x).$$

Here, $\sup |\xi_{\alpha(i)}(u, x)| \leq$ const., the $c_i$'s are constants and $\alpha(i), \beta(i) \in A$ with $|\alpha(i)| \geq \ell$ and $|\beta(i)| \geq k - \ell + 1$. Assuming, for simplicity, that the design points are uniformly distributed, we may construct the following root-$n$ consistent estimator of $g_1(x)$ using the approach at (2.3):

$$\hat{g}_1(x) = \sum_{i=1}^{r} \frac{1}{n_{\alpha(i)}} \sum_{j=1}^{n_{\alpha(i)}} Y_{\alpha(i)j} \xi_{\alpha(i)}(X_{\alpha(i)j}, x).$$

We may estimate $g_2(x)$, using the method at (2.8), as follows:

$$\hat{g}_2(x) = \sum_{i=1}^{s} \frac{c_i}{n_{\alpha_0} h^{k - |\beta(i)|}} \sum_{j=1}^{n_{\alpha_0}} Y_{\alpha_0 j} \, K\left(\frac{X^*_{\beta(i)j} - x^*}{h}\right).$$

Here $X^*_{\beta(i)j}$ and $x^*$ denote the vectors of those $k - |\beta(i)|$ components of $X_{\beta(i)j}$ and $x$, respectively, for which the corresponding components of $\beta(i)$ are zero. Since $k - |\beta(i)| \leq \ell - 1$ for each $i$, the squared-error convergence rate of $\hat{g}_2$ to $g_2$ is $n^{-2d/(2d+\ell-1)}$. Therefore, the squared-error convergence rate of $\hat{g} = \hat{g}_1 + \hat{g}_2$ to $g$ is also $n^{-2d/(2d+\ell-1)}$, as claimed three paragraphs above.

EXAMPLE 5. Suppose $k = 2$ and that noisy data are available for $g^{(1,1)}$ and $g^{(0,0)} = g$. Use the root-$n$ consistent estimator of $\underline{g}_{11}$ from Example 2 to write

$$y_{(0,0)} - \hat{\underline{g}}_{11}(x_1, x_2) = g(x_1, 0) + g(0, x_2) - g(0,0) + O_p(n^{-1/2}) + \varepsilon_{(0,0)},$$

which is additively separable in $x_1$ and $x_2$ and hence estimable at one-dimensional optimal rates.

EXAMPLE 6. Suppose $k = 3$ and that noisy data are available for $g^{(1,1,1)}$, $g^{(1,1,0)}$, $g^{(1,0,1)}$, $g^{(0,1,1)}$ and $g^{(0,0,0)} = g$. Define

$$\underline{g}_{111}(x_1, x_2, x_3) \equiv \int_0^{x_1}\int_0^{x_2}\int_0^{x_3} g^{(1,1,1)}(u_1, u_2, u_3)\, du_1\, du_2\, du_3$$
$$= g(x_1, x_2, x_3) - g(x_1, x_2, 0) - g(x_1, 0, x_3) - g(0, x_2, x_3)$$
$$+ g(x_1, 0, 0) + g(0, x_2, 0) + g(0, 0, x_3) - g(0, 0, 0),$$

$$\underline{g}_{110}(x_1, x_2, \cdot) \equiv \int_0^{x_1}\int_0^{x_2}\int_0^1 g^{(1,1,0)}(u_1, u_2, x_3)\, du_1\, du_2\, dx_3$$
$$- \int_0^1 \underline{g}_{111}(x_1, x_2, x_3)\, dx_3$$
$$= g(x_1, x_2, 0) - g(x_1, 0, 0) - g(0, x_2, 0) + g(0, 0, 0),$$

$$\underline{g}_{101}(x_1, \cdot, x_3) \equiv \int_0^{x_1}\int_0^1\int_0^{x_3} g^{(1,0,1)}(u_1, x_2, u_3)\, du_1\, dx_2\, du_3$$
$$- \int_0^1 \underline{g}_{111}(x_1, x_2, x_3)\, dx_2$$
$$= g(x_1, 0, x_3) - g(x_1, 0, 0) - g(0, 0, x_3) + g(0, 0, 0),$$

$$\underline{g}_{011}(\cdot, x_2, x_3) \equiv \int_0^1\int_0^{x_2}\int_0^{x_3} g^{(0,1,1)}(x_1, u_2, u_3)\, dx_1\, du_2\, du_3$$
$$- \int_0^1 \underline{g}_{111}(x_1, x_2, x_3)\, dx_1$$
$$= g(0, x_2, x_3) - g(0, x_2, 0) - g(0, 0, x_3) + g(0, 0, 0).$$

Sample analogues of all integral expressions may be calculated without local averaging. Thus, $\underline{g}_{111}$ and its integrals can be estimated root-$n$ consistently, as can $\underline{g}_{110}, \underline{g}_{101}$ and $\underline{g}_{011}$. Now, write

$$y_{(0,0,0)} - \hat{\underline{g}}_{111}(x_1, x_2, x_3) - \hat{\underline{g}}_{110}(x_1, x_2, \cdot) - \hat{\underline{g}}_{101}(x_1, \cdot, x_3) - \hat{\underline{g}}_{011}(\cdot, x_2, x_3)$$
$$= g(x_1, 0, 0) + g(0, x_2, 0) + g(0, 0, x_3) - 2g(0, 0, 0) + O_p(n^{-1/2}) + \varepsilon_{(0,0,0)},$$

which is additively separable in $x_1, x_2$ and $x_3$ and hence estimable at one-dimensional optimal rates.

### 3. Estimation.

3.1. *Smoothing techniques.* In Section 2, we gave examples of estimators in the case where the design points $X_{\alpha i}$ are uniformly distributed on $\mathcal{R}_k$. More generally, we should normalize the summands of our estimators, such as those at (2.3) and (2.8), using estimators of the densities of the distributions of design points. For simplicity, we shall develop the case of (2.8) in this setting, noting that other cases are similar.

Suppose we observe the datasets at (2.7) for each $\beta \in B$, where the latter is the set of $p$-vectors of zeros and ones with $1 \leq p \leq k$. Note that we may also interpret $\beta$ as a $k$-vector, an element of $A$, in which each of the last $k - p$ components is zero. Both interpretations will be made below.

The design points $X_{\beta i}$, which are $k$-vectors, are assumed to be distributed on $\mathcal{R}_k$ with density $f_\beta$, say. As in Section 2.2, let $X_{\beta i}^{[k-p]}$ and $x^{[k-p]}$ denote the $(k - p)$-vectors consisting of the last $k - p$ components of $X_{\beta i}$ and $x$, respectively, let $K$ be a $(k - p)$-variate kernel function, let $h$ denote a bandwidth and redefine

$$(3.1) \quad (\widehat{\psi_\beta g^\beta})(x) = \frac{1}{n_\beta h^{k-p}} \sum_{i=1}^{n_\beta} \frac{Y_{\beta i} \chi_\beta(X_{\beta i}, x)}{\tilde{f}_{\beta, -i}(X_{\beta i})} K\left(\frac{X_{\beta i}^{[k-p]} - x^{[k-p]}}{h}\right),$$

where $\tilde{f}_{\beta, -i}$ denotes an estimator of $f_\beta$ computed from the dataset $\mathcal{X}_{\beta, -i} = \{X_{\beta 1}, \ldots, X_{\beta n_\beta}\} \setminus \{X_{\beta i}\}$ obtained by dropping the $i$th observation. Note that $\chi_\beta(X_{\beta i}, x)$ depends only on the first $p$ components of $X_{\beta i}$ and $x$, whereas $\tilde{f}_{\beta, -i}(x)$ and $f_\beta(x)$ depend nondegenerately on all $k$ components of $x$.

A degree of interest centers on the definition adopted for $\tilde{f}_{\beta, -i}$. We shall discuss an edge-corrected kernel method, but, of course, there are many other techniques that can be used—for example, polynomial interpolation, or polynomial smoothing, applied to binned data.

Let $H > 0$ denote a bandwidth and let $L_1$ represent a bounded function of a real variable $t$, supported on the interval $[-1, 1]$ and satisfying $\int t^j K_1(t)\, dt = \delta_{j0}$ (the Kronecker delta), for $0 \leq j \leq d_1 - 1$. (The positive integer $d_1$ may differ from the order $d$ of the kernel $K$.) Construct a $k$-variate product kernel $L$,

$$(3.2) \qquad L(u_1, \ldots, u_k) = L_1(u_1) \cdots L_1(u_k).$$

The density estimator

$$(3.3) \qquad \hat{f}_{\beta, -i}(x) = \frac{1}{n_\beta - 1} \sum_{j : j \neq i} L\left(\frac{x - X_{\beta j}}{H}\right)$$

does not suffer edge effects provided $x_i \in [h, 1 - h]$ for $1 \leq i \leq k$. However, if for one or more values of $i$, $x_i$ lies outside $[h, 1 - h]$ and, more particularly, if $0 < x_i < h$, then edge effects may be averted by replacing $L_1(u_i)$ with $L_{\text{edge}}(u_i)$ in the definition of $L$ at (3.2). Here, $L_{\text{edge}}$ is a bounded, univariate edge kernel, supported on $[0,1]$ and satisfying $\int t^j L_{\text{edge}}(t)\,dt = \delta_{j0}$ for $0 \leq j \leq d_1 - 1$. Similarly, if $1 - h < x_i < 1$, then we would replace $L_1$, applied to the $i$th component in (3.2), by an edge kernel supported on $[-1, 0]$.

With these modifications, the density estimator $\hat{f}_{\beta,-i}$ defined at (3.3) is of $d_1$th order and does not suffer edge effects in $\mathcal{R}_k$.

Our definition of $\tilde{f}_{\beta,-i}$ ensures that the estimator at (3.1) is protected from edge effects in the first $p$ coordinates of $x$. However, we should modify $K$ in the same way as we did $L$; otherwise, $\widehat{\psi_\beta g^\beta}$ will suffer from edge effects in the last $k-p$ coordinates of $x$. We shall assume that this has been done so that the $(k-p)$-variate kernel $K$ is, analogously to $L$, a product of $k-p$ bounded, compactly supported, $d$th order univariate kernels that are switched to appropriate edge kernels if one or more components of $x^{[k-p]}$ are within $h$ of the boundary. The univariate kernels, $K_1$ and $K_{\text{edge}}$, say, will each be taken to be of $d$th order.

Rather than employ special kernels to overcome edge effects, we may use local polynomial methods to construct $\widehat{\psi_\beta g^\beta}$, obtaining an alternative estimator to that at (3.1). In this approach, we would run a $(k-p)$-variate local polynomial smoother of degree $d-1$ through the data pairs

$$(3.4) \qquad (X_{\beta i}^{[k-p]}, Y_{\beta i} \chi_\beta(X_{\beta i})/\tilde{f}_{\beta,-i}(X_{\beta i})), \qquad 1 \leq i \leq n_\beta.$$

This technique is also able to correct for a nonuniform joint distribution of the last $k-p$ components, so we could normalize the "response variable" a little differently than by dividing by $\tilde{f}_{\beta,-i}(X_{\beta i})$, as at (3.4). However, the normalization at (3.4) causes no problems for the local polynomial smoother.

3.2. *Limit theory for estimators.* For the sake of simplicity, we shall give theory only for edge-corrected kernel approaches to estimation. In particular, we assume $\tilde{f}_{\beta,-i}$ is constructed using the methods described in Section 3.1, that the univariate kernel $L_1$ and its two edge-correcting forms $L_{\text{edge}}$ are bounded and compactly supported and that the same is true of the univariate kernels $K_1$ and $K_{\text{edge}}$ that are multiplied together to give the $(k-p)$-variate kernel $K$. To this, we add the assumption that

(3.5) $\quad K_1, K_{\text{edge}}, L_1$ and $L_{\text{edge}}$ are Hölder

continuous as functions on the real line.

Recall that the estimator $\tilde{f}_{\beta,-i}$ is constructed using a $d_1$th order kernel $L$ and a bandwidth $H$ and that the kernel $K$ used in (3.1) is of order $d$ and employs a bandwidth $h$. Of these quantities, we assume the following conditions:

$$(3.6) \qquad d > \tfrac{1}{2}(k+p) \quad \text{and} \quad d_1 > k,$$

(3.7) for constants $0 < C_1 < C_2 < \infty$ and $\eta > 0$,
$C_1 n_\beta^{-1/(2d+k-p)} \leq h \leq C_2 n_\beta^{-1/(2d+k-p)}$ and
$C_1 n_\beta^{\{-1/(2k)\}+\eta} \leq H \leq C_2 \min\{n_\beta^{-1/(2d_1)}, n_\beta^{-1/(2d+k-p)}\} n_\beta^{-\eta}$
for all sufficiently large $n_\beta$.

Provided (3.6) holds, we may choose $h$ and $H$ satisfying (3.7). We also suppose that

(3.8) $g^\beta$ is bounded, the last $k - p$ components of $g$ have $d$ continuous derivatives and $f_\beta$ has $d_1$ bounded derivatives and is bounded away from zero on $\mathcal{R}_k$.

We also make the following basic "structural" assumptions:

(3.9) data pairs $(X_{\beta i}, Y_{\beta i})$ are generated by the model at (2.7), in which the design variables $X_{\beta i}$ are independent and identically distributed on $\mathcal{R}_k$ with density $f_\beta$, the errors $\varepsilon_{\beta i}$ are independent and identically distributed with zero mean and the errors are independent of the design points.

From these data, construct the estimator $\widehat{\psi_\beta g^\beta}$ defined at (3.1). Recall that $u^{[k-p]}$ denotes the $(k-p)$-vector consisting of the last $k-p$ components of the $k$-vector $u$. Let $w(u, x \mid h)$ represent the $k$-vector with $u_j$ in position $j$ for $1 \leq j \leq p$ and $x_j + h_j u_j$ in position $j$ for $p+1 \leq j \leq k$.

THEOREM 3. *Assume* $1 \leq p \leq k$, *that conditions* (3.5)–(3.9) *hold and that the distribution of the errors* $\varepsilon_{\beta i}$ *has zero mean and all moments finite. Then*

$$(3.10) \quad \begin{aligned} (\widehat{\psi_\beta g^\beta})(x) &= \int_{\mathcal{R}_k} g\{w(u, x \mid h)\} \chi_\beta\{w(u, x \mid h), x\} K(u^{[k-p]}) \, du \\ &\quad + \frac{1}{n_\beta} \sum_{i=1}^{n_\beta} \frac{\varepsilon_{\beta i} \chi_\beta(X_{\beta i}, x)}{f_\beta(X_{\beta i}) h^{k-p}} K\left(\frac{X_{\beta i}^{[k-p]} - x^{[k-p]}}{h}\right) \\ &\quad + o_p(n_\beta^{-d/(2d+k-p)}), \end{aligned}$$

*uniformly in* $x \in \mathcal{R}_k$.

We shall discuss the implications of Theorem 3 in the two main cases, $p = k$ and $p < k$. In the first setting, the contribution of the kernel $K$ to (3.10) is degenerate and the integral on the right-hand side is identical to $(\psi_\beta g^\beta)(x)$. (Here, $\beta$ is a $k$-vector.) Therefore, when $p = k$, (3.10) is equivalent to

$$(3.11) \quad (\widehat{\psi_\beta g^\beta})(x) = (\psi_\beta g^\beta)(x) + Z_{n_\beta}(x) + o_p(n_\beta^{-1/2}),$$

uniformly in $x \in \mathcal{R}_k$, where

$$Z_{n_\beta}(x) = \frac{1}{n_\beta} \sum_{i=1}^{n_\beta} \frac{\varepsilon_{\beta i} \chi_\beta(X_{\beta i}, x)}{f_\beta(X_{\beta i})}$$

is a zero-mean stochastic process defined on $\mathcal{R}_k$. As $n_\beta$ increases, $n_\beta^{1/2} Z_{n_\beta}$ converges weakly to the Gaussian process $Z_0$, say, with zero mean and covariance function

$$(3.12) \quad \text{cov}\{Z_0(x_1), Z_0(x_2)\} = \sigma_\beta^2 \int_{\mathcal{R}_k} \chi_\beta(u, x_1) \chi_\beta(u, x_2) f_\beta(u)^{-1} \, du,$$

where $\sigma_\beta^2 = \text{var}(\varepsilon_{\beta i})$. This property and (3.11) together imply that $\widehat{\psi_\beta g^\beta}$ converges uniformly to $\psi_\beta g^\beta$ at rate $n^{-1/2}$:

$$\sup_{x \in \mathcal{R}_k} |(\widehat{\psi_\beta g^\beta})(x) - (\psi_\beta g^\beta)(x)| = O_p(n_\beta^{-1/2}).$$

Next, we treat the case $p < k$. Although $\chi_\beta(u, x)$ is discontinuous as a function of the first $p$ components of $u$, if $g$ has $d$ continuous derivatives of its last $k - p$ components, then so too does $\chi_\beta(\cdot, x)$; see the definition of $\chi_\alpha$ given in Appendix A.1 and recall that definition has a minor adaptation to the case of $\chi_\beta$. Therefore, standard Taylor expansion methods may be used to prove that for a continuous function $a$,

$$(3.13) \quad \begin{aligned} &\int g\{w(u, x \,|\, h)\} \chi_\beta\{w(u, x \,|\, h), x\} K(u^{[k-p]}) \, du \\ &= \int g\{w(u, x), x\} \chi_\beta\{w(u, x), x\} \, du_1 \cdots du_p + h^d \, a(x) + o(h^d) \end{aligned}$$

as $h \to 0$, where $w(u, x) = w(u, x \,|\, 0)$ is the $k$-vector with $u_j$ in position $j$ for $1 \leq j \leq k$ and $x_j$ in position $j$ for $p + 1 \leq j \leq k$. The series on the right-hand side of (3.10) is asymptotically normally distributed with zero mean and variance $(n_\beta h^{k-p})^{-1} \sigma_\beta^2 \tau(x)^2 \kappa$, where

$$\tau(x)^2 = \int \chi_\beta\{w(u, x), x\}^2 f_\beta\{w(u, x)\}^{-1} \, du_1 \cdots du_p$$

and $\kappa = \int K^2$. This result, (3.10) and (3.13) collectively imply that for the choice of $h$ given in (3.7), $(\widehat{\psi_\beta g^\beta})(x)$ converges to $(\psi_\beta g_\beta)(x)$ at the pointwise squared-error rate $n_\beta^{-2d/(2d+k-p)}$, as claimed in Section 2.2. The uniform convergence rate is slower only by a logarithmic factor.

It is straightforward to prove that the pointwise rate is minimax optimal. Indeed, that property follows from conventional minimaxity results in nonparametric regression on taking $g$ to be a function of which the dependence on the first $p$ coordinates is degenerate. Likewise, the uniform convergence rate can be shown to be optimal, provided we use a slightly larger bandwidth $h$, increased by a logarithmic factor relative to that asserted in (3.7).

We close by formally stating the main results discussed above.

COROLLARY 2. *Assume the conditions of Theorem 3. If $p = k$, then $n_\beta^{1/2}\{\widehat{\psi_\beta g^\beta}(x) - \psi_\beta g^\beta(x)\}$, viewed as a stochastic process indexed by $x \in \mathcal{R}_k$, converges weakly, as $n_\beta \to \infty$, to a zero-mean Gaussian process $Z_0$ with covariance function given at (3.12). If $p < k$ and if $h \sim \text{const} \cdot n^{-1/(2d+k-p)}$, then for each $x \in \mathcal{R}_k$, $n_\beta^{d/(2d+k-p)}\{\widehat{\psi_\beta g^\beta}(x) - \psi_\beta g^\beta(x)\}$ is asymptotically normally distributed with finite mean and variance.*

Of course, in order to construct an estimator $\hat{g}$ of $g$, we must add $\widehat{\psi_\beta g^\beta}$ over all $\beta$; see (2.9). The resulting limit theory for $\hat{g}$ is the superposition of that for each $\widehat{\psi_\beta g^\beta}$. However, provided the sets of design points $X_{\beta i}$ and errors $\varepsilon_{\beta i}$ are independent for different $\beta$'s, properties of the superposition are readily derived from the results that we have already obtained for a single $\beta$.

Indeed, under this assumption of row-wise independence, it follows directly from Corollary 2 that if, for a sequence of integers $n$ diverging to infinity, $n_\beta/n$ converges to a strictly positive constant $c_\beta$ for each $\beta \in B$, then (a) if $p = k$, $n^{1/2}(\hat{g} - g)$ converges weakly to a zero-mean Gaussian process defined on $\mathcal{R}_k$ and (b) if $p < k$, then for each $x \in \mathcal{R}_k$, $n^{d/(2d+k-p)}\{\hat{g}(x) - g(x)\}$ is asymptotically normally distributed with finite mean and variance.

Correlation among residuals in different equations can also be accommodated. Let $B = \{\beta_1, \ldots, \beta_s\}$. Suppose $(X_i, Y_{\beta_1 i}, \ldots, Y_{\beta_s i})_{i=1,\ldots,n}$ are independent and identically distributed, where $Y_{\beta_j i} = g^{\beta_j}(X_i) + \varepsilon_{\beta_j i}$, $j = 1, \ldots, s$ and $\sigma_{jj'} \equiv \text{cov}(\varepsilon_{\beta_j i}, \varepsilon_{\beta_{j'} i})$. Let $f(x)$ denote the design density of the $X_i$ which are distributed independently of the residuals. Then conclusions (a) and (b) of the previous paragraph continue to hold with the covariance function of the limiting Gaussian process in (a), say $Z_0$, given by

$$\text{cov}\{Z_0(x_1), Z_0(x_2)\} = \sum_{jj'} \sigma_{jj'} \int_{\mathcal{R}_k} \chi_{\beta_j}(u, x_1) \chi_{\beta_{j'}}(u, x_2) f(u)^{-1} du.$$

## 4. Numerical results.

4.1. *Simulation of cost function and input factor estimation.* We return to the cost function estimation problem discussed in Section 1. Since doubling of input prices at a given level of output doubles costs, the cost function is homogeneous of degree one in input prices. Thus, we may write average costs, that is, costs per unit of output $Q$, as $\text{AC} = r\,g(Q, w/r)$, where $r$ and $w$ are the prices of capital and labor, respectively. Applying Shephard's Lemma yields the average labor function, $\mathrm{AL} = \partial \mathrm{AC}/\partial w = \partial g(Q, w/r)/\partial(w/r)$. If noisy data are available for AC and AL, then this application is analogous to Example 3 above, except that the nonparametric function $g$ is multiplied by $r$, a feature which arises from the degree-one homogeneity of cost functions in their factor prices.

We calibrate our simulations using the Cobb–Douglas production function $Q = cK^{c_1}L^{c_2}$ (see, e.g., [23]). The data-generating mechanism for average costs is

$$(4.1) \qquad y_{(0,0)} = \mathrm{AC} + \varepsilon_{(0,0)} = r\tilde{c}Q^{\frac{1-c_1-c_2}{c_1+c_1}}\left(\frac{w}{r}\right)^{\frac{c_2}{c_1+c_2}} + \varepsilon_{(0,0)},$$

where $\tilde{c} = ((c_1/c_2)^{c_2/(c_1+c_2)} + (c_1/c_2)^{-c_1/(c_1+c_2)})c^{\frac{-1}{c_1+c_2}}$. For average labor, we use

$$(4.2) \qquad y_{(0,1)} = \mathrm{AL} + \varepsilon_{(0,1)} = \frac{c_2}{c_1+c_2}\tilde{c}Q^{\frac{1-c_1-c_2}{c_1+c_2}}\left(\frac{w}{r}\right)^{\frac{-c_1}{c_1+c_2}} + \varepsilon_{(0,1)}.$$

In the simulations below, we set $c_1 = 0.8$ and $c_2 = 0.7$. Data for $Q$ and for the ratio of factor prices $w/r$ are generated from independent uniform distributions on $[0.5, 1.5]$. We assume that $\varepsilon_{(0,0)}$ and $\varepsilon_{(0,1)}$ are normal residuals with zero means, standard deviations 0.35 and correlation $\rho$ set to 0.0, 0.4 or 0.9. The $R^2$ is approximately 0.75 for the AC equation and 0.15 for the AL equation. Our reference estimator of average costs consists of applying bivariate kernel smoothing to the triples $(y_{(0,0)}/r, Q, w/r)$ to obtain $\hat{g}(Q, w/r)$, which is then multiplied by $r$.

To incorporate the labor data, define

$$(4.3) \qquad \hat{\underline{g}}_a(r, Q, w/r) = r\frac{1}{nh} \sum_{\substack{|Q_j-Q|\leq h/2 \\ w_j/r_j \leq w/r}} Y_{(0,1)j},$$

$$(4.4) \qquad \hat{\underline{g}}_b(r, Q, \cdot) = r\frac{1}{nh} \sum_{|Q_j-Q|\leq h/2} \frac{Y_{(0,0)j} - \hat{\underline{g}}_a(r_j, Q_j, w_j/r_j)}{r_j}.$$

Then $\widehat{\mathrm{AC}} = \hat{\underline{g}}_a + \hat{\underline{g}}_b$. Table 1 summarizes our results for various sample sizes $n$ and residual correlations. There, we report the mean squared errors of this estimator relative to the bivariate kernel estimator described above. There are substantial efficiency gains, which increase with sample size, as would be expected given the faster convergence rates of derivative-based estimators.

4.2. *Estimating costs of electricity distribution.* To further illustrate the procedure, we use data on 81 electricity distributors in Ontario. (For additional details, see [24].) We have data on output, $Q$, which is the number of customers served and which varies from about 500 to over 200,000. Average labor, AL, equals the number of employees divided by $Q$. In addition, we have data on hourly wages, $w$, and the cost of capital, $r$.

TABLE 1
*MSEs of derivative-based AC estimator relative to bivariate kernel smoothing*

| $n$ | $\rho = 0.0$ | $\rho = 0.4$ | $\rho = 0.9$ |
|------|--------------|--------------|--------------|
| 100  | 0.384        | 0.384        | 0.374        |
| 200  | 0.277        | 0.274        | 0.272        |
| 500  | 0.233        | 0.230        | 0.228        |
| 1000 | 0.185        | 0.187        | 0.186        |

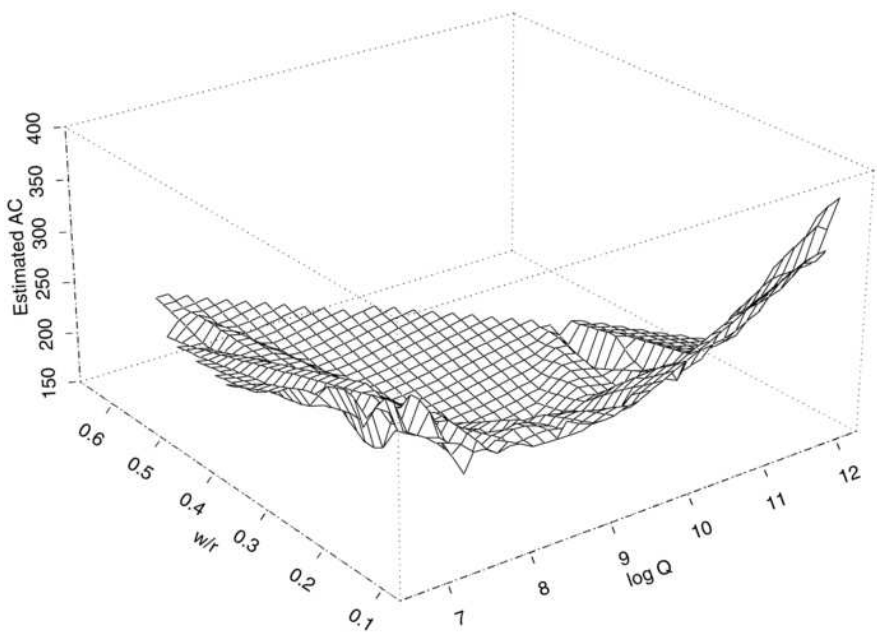

FIG. 1. *Function estimate using data on function only.*

Figure 1 illustrates the estimated average cost function using only AC data and a bivariate loess smoother available in S-PLUS. Next, we use both the AC and AL data and apply equations (4.3) and (4.4), suitably modified for the nonuniform distribution of $w/r$. Figure 2 illustrates the resulting estimate.

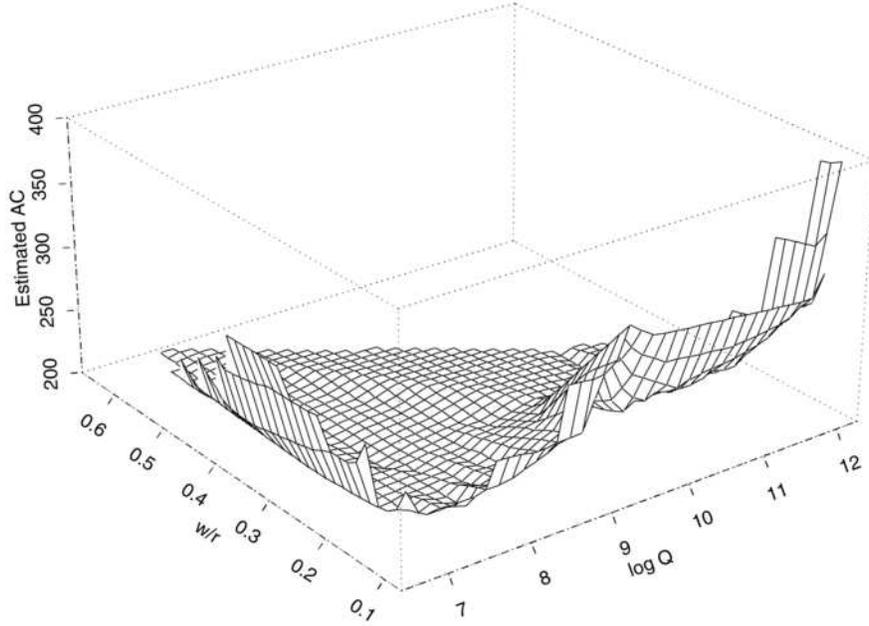

Fig. 2. *Function estimate using data on function and partial derivative.*

## APPENDIX: TECHNICAL ARGUMENT

**A.1. Proof of Theorem 1.** It is readily seen that when $k = 1$,

$$(A.1) \qquad g(x) = \sum_{j=0}^{1} \int_0^1 \chi_j(u,x), g^{(j)}(u)\, du,$$

where $\chi_0(u,x) \equiv 1$, $\chi_1(u,x) \equiv u - 1 + I(u \le x)$, $I(u \le x) = 1$ if $u \le x$ and equals 0 otherwise and $g^{(j)}(x) = (\partial/\partial x)^j g(x)$. Repeating identity (A.1) for each component of a function $g$ of $k \ge 1$ variables, we deduce that Theorem 1 holds with $\psi_\alpha$ defined by $(\psi_\alpha g)(x) = \int_{\mathcal{R}_k} \chi_\alpha(u,x) g(u)\, du$, where

$$\chi_\alpha(u_1,\ldots,u_k, x_1,\ldots,x_k) = \prod_{j=1}^{k} \chi_{\alpha_j}(u_j, x_j)$$

and $\alpha = (\alpha_1,\ldots,\alpha_k)$. Note, particularly, that $|\chi_\alpha| \le 1$ and so $\psi_\alpha \in \mathcal{K}(1)$, where $\mathcal{K}(C)$ is defined as in Section 2.1.

**A.2. Proof of Theorem 2.** In proving the theorem, we may assume that $\mathcal{P} = \{1,\ldots,k\}$, since the contrary case can be treated by fixing components of which the index does not lie in $\mathcal{P}$. In the notation at (2.11), define

$$(N_\alpha b)(x) = \int_0^{x_{i_1}} \cdots \int_0^{x_{i_{|\alpha|}}} b\{v_\alpha(u,x)\}\, du_{i_1} \cdots du_{i_{|\alpha|}}.$$

Given $\alpha \in A$, let $A(\alpha)$ denote the set of vectors $\beta = (\beta_1, \ldots, \beta_k) \in A$ for which each index $j$ with $\beta_j = 1$ is also an index with $\alpha_j = 1$. Put $\alpha_0 = (0, \ldots, 0)$ and $A_1(\alpha) = A(\alpha) \setminus \{\alpha_0\}$. We shall prove shortly that for all $\alpha \in A$ and $b \in \mathcal{B}_k$,

$$(A.2) \qquad \sum_{\beta \in A(\alpha)} (-1)^{|\beta|} M_\beta N_\alpha b^\alpha = \sum_{\beta \in A(\alpha)} (-1)^{|\beta|} M_\beta b$$

or, equivalently,

$$(A.3) \qquad b = \sum_{\beta \in A(\alpha)} (-1)^{|\beta|} M_\beta N_\alpha b^\alpha - \sum_{\beta \in A_1(\alpha)} (-1)^{|\beta|} M_\beta b.$$

Substituting $b = g$ and $\alpha = (1, \ldots, 1)$ into (A.3), we obtain

$$(A.4) \qquad g = \sum_{\beta \in A} (-1)^{|\beta|} M_\beta N_\alpha g^\alpha - \sum_{\beta \in A_1(\alpha)} (-1)^{|\beta|} M_\beta g,$$

where $A_1 = A_1(1, \ldots, 1) = A \setminus \{\alpha_0\}$. The first series on the right-hand side is a linear expression in integrals of the form at (2.12). If $|\beta| \geq k - \ell + 1$, then $M_\beta g$ is also of the form claimed in the theorem. It remains only to treat terms $M_\beta g$ with $|\beta| \leq k - \ell$, which we do using an iterative argument. [Note that, since $\beta \in A_1(\alpha)$, we have $|\beta| \geq 1$, so we have already finished if $\ell = k$.]

Write $\mathcal{S}(\beta)$ for the set of indices $i$ such that $\beta_i = 1$ and define $\alpha^1 = \alpha^1(\beta) \in A$ by $\mathcal{S}(\alpha_1) = \mathcal{S}(\alpha) \setminus \mathcal{S}(\beta)$. Apply (A.3) again, this time with $\alpha = \alpha^1$ and $b = M_\beta g$, obtaining

$$M_\beta g = \sum_{\beta^1 \in A(\alpha^1)} (-1)^{|\beta^1|} M_{\beta^1} N_{\alpha^1} (M_\beta g)^{\alpha^1} - \sum_{\beta^1 \in A_1(\alpha^1)} (-1)^{|\beta^1|} M_{\beta^1} M_\beta g.$$

By definition of $\alpha^1$, $(M\beta g)^{\alpha^1} = (M\beta)(g^{\alpha^1})$, and so $N_{\alpha^1}(M_\beta g)^{\alpha^1} = (N_{\alpha^1} M_\beta)(g^{\alpha^1})$, which is a $k$-fold integral of $g^{\alpha^1}$, where $|\alpha^1| = k - |\beta| \geq k - (k - \ell) = \ell$. Also, $M_{\beta^1} M_\beta = M_{\beta^2} g$, where $\beta^2 \in A$ and $|\beta^2| \geq 2$. (The superscript 2 is an index, not an exponent.) If $|\beta^2| \geq k - \ell + 1$, we are done; if $|\beta^2| \leq k - \ell$, we continue the process of iteration.

Finally, we derive (A.2). Again, it suffices to treat the case $\alpha = (1, \ldots, 1)$, since other contexts may be addressed by fixing components $x_j$ for $j$ such that $\alpha_j = 0$. In the case $\alpha = (1, \ldots, 1)$,

$$(N_\alpha b^\alpha)(x) = \int_0^{x_1} \cdots \int_0^{x_k} b^{(1,\ldots,1)}(u_1, \ldots, u_k) \, du_1 \cdots du_k$$

$$= \sum_{\gamma \in A} (-1)^{|\gamma|} b\{v_\gamma(0, x)\},$$

whence

$$(A.5) \quad \sum_{\beta \in A} (-1)^{|\beta|} (M_\beta N_\alpha b^\alpha)(x) = \sum_{\gamma \in A} (-1)^{|\gamma|} \sum_{\beta \in A} (-1)^{|\beta|} [M_\beta b\{v_\gamma(0, \cdot)\}](x).$$

If $\gamma_i = 1$ and $\beta = (\beta_1, \ldots, \beta_k) \in A$ then, if we switch $\beta_i$ from 0 to 1, we do not alter the value of $[M_\beta b\{v_\gamma(0, \cdot)\}](x)$. Therefore, by virtue of the factor $(-1)^{|\beta|}$ below,

$$\sum_{\beta \in A} (-1)^{|\beta|} [M_\beta b\{v_\gamma(0, \cdot)\}] \equiv 0$$

unless $\gamma = \alpha_0$. However, $v_{\alpha_0}(0, u) = u$, so by (A.5),

$$\sum_{\beta \in A} (-1)^{|\beta|} (M_\beta N_\alpha b^\alpha)(x) = \sum_{\beta \in A} (-1)^{|\beta|} (M_\beta b)(x),$$

which, in the case $\alpha = (1, \ldots, 1)$, is equivalent to (A.2).

**A.3. PROOF OF THEOREM 3.** The estimators $\hat{f}_{\beta, -i}$ have biases and variances that are uniformly of orders $H^{d_1}$ and $(n_\beta H^k)^{-1}$, respectively, and, in particular,

(A.6) $$\sup_{x \in \mathcal{R}_k, 1 \leq i \leq n_\beta} |E\{\hat{f}_{\beta,-i}(x)\} - f(x)| = O(H^{d_1}).$$

Arguments based on Markov's inequality show that for each $c, C > 0$,

(A.7) $$\sup_{x \in \mathcal{R}_k, 1 \leq i \leq n_\beta} P\{|\hat{f}_{\beta,-i}(x) - E\hat{f}_{\beta,-i}(x)| > (n_\beta^{c-1} H^{-k})^{1/2}\} = O(n_\beta^{-C}).$$

The Hölder continuity assumed of $L$ may be used to prove that if $C_1 > 0$ is chosen sufficiently large, then for all $C_2 > 0$,

$$E\left\{\sup_{|x_1 - x_2| \leq n_\beta^{-C_1}, 1 \leq i \leq n_\beta} |\hat{f}_{\beta,-i}(x_1) - \hat{f}_{\beta,-i}(x_2)|^{C_2}\right\} = O(n_\beta^{-C_2}).$$

Therefore, again by Markov's inequality and for each $c, C > 0$,

(A.8) $$P\left\{\sup_{|x_1 - x_2| \leq n_\beta^{-C_1}, 1 \leq i \leq n_\beta} |\hat{f}_{\beta,-i}(x_1) - \hat{f}_{\beta,-i}(x_2)| > n^{c-1}\right\} = O(n_\beta^{-C}).$$

Applying (A.7) on a lattice of values $x \in \mathcal{R}_k$ of edge width $n_\beta^{-C_1}$ and using (A.8) to bound $|\hat{f}_{\beta,-i}(x_1) - \hat{f}_{\beta,-i}(x_2)|$ when $x_1$ is off the lattice and $x_2$ is the nearest grid point to $x_1$, we may prove that for each $c, C > 0$,

(A.9) $$P\left\{\sup_{x \in \mathcal{R}_k, 1 \leq i \leq n_\beta} |\hat{f}_{\beta,-i}(x) - E\hat{f}_{\beta,-i}(x)| (n_\beta^{c-1} H^{-k})^{1/2}\right\} = O(n_\beta^{-C}).$$

Below, we shall refer to this as the "lattice argument"; it employs the Hölder-continuity condition (3.5).

Taylor expanding $\hat{f}_{\beta,-i}^{-1}$ as $\hat{f}_{\beta,-i}^{-1} = f_\beta^{-1} - (\hat{f}_{\beta,-i} - f_\beta) f_\beta^{-2} + \cdots$, we may show that

(A.10) $$\tilde{f}_{\beta,-i}(X_{\beta i})^{-1} - f_\beta(X_{\beta i})^{-1} = -\frac{\hat{f}_{\beta,-i}(X_{\beta i}) - f_\beta(X_{\beta i})}{f_\beta(X_{\beta i})^2} + \Delta_{\beta i},$$

where, by (A.6) and (A.9), we have for each $c > 0$,

(A.11) $$\max_{1 \leq i \leq n_\beta} |\Delta_{\beta i}| = O_p(n_\beta^{c-1} H^{-k} + H^{d_1}).$$

Substituting (A.10) into the definition (3.1) of the estimator $\widehat{\psi_\beta g^\beta}$, we deduce that

(A.12) $$(\widehat{\psi_\beta g^\beta})(x) = S_1(x) - S_2(x) - S_3(x) - S_4(x) + S_5(x),$$

where

$$S_1(x) = \frac{1}{n_\beta} \sum_{i=1}^{n_\beta} \frac{Y_{\beta i} \chi_\beta(X_{\beta i}, x)}{f_\beta(X_{\beta i})} K_{\beta i}(x),$$

$$S_2(x) = \frac{1}{n_\beta} \sum_{i=1}^{n_\beta} \frac{g(X_{\beta i}) \chi_\beta(X_{\beta i}, x) \{\hat{f}_{\beta,-i}(X_{\beta i}) - \kappa_\beta(X_{\beta i})\}}{f_\beta(X_{\beta i})^2} K_{\beta i}(x),$$

$$S_3(x) = \frac{1}{n_\beta} \sum_{i=1}^{n_\beta} \frac{\varepsilon_{\beta i} \chi_\beta(X_{\beta i}, x) \{\hat{f}_{\beta,-i}(X_{\beta i}) - \kappa_\beta(X_{\beta i})\}}{f_\beta(X_{\beta i})^2} K_{\beta i}(x),$$

$$S_4(x) = \frac{1}{n_\beta} \sum_{i=1}^{n_\beta} \frac{Y_{\beta i} \chi_\beta(X_{\beta i}, x) \{\kappa_\beta(X_{\beta i}) - f_\beta(X_{\beta i})\}}{f_\beta(X_{\beta i})^2} K_{\beta i}(x),$$

$$S_5(x) = \frac{1}{n_\beta} \sum_{i=1}^{n_\beta} Y_{\beta i} \chi_\beta(X_{\beta i}, x) \Delta_{\beta i} K_{\beta i}(x),$$

$K_{\beta i}(x) = h^{-(k-p)} K\{(X_{\beta i}^{[k-p]} - x^{[k-p]})/h\}$ and $\kappa_\beta(x) = E\{\hat{f}_{\beta,-i}(x)\}$.

Noting that the errors $\varepsilon_{\beta i}$ are independent of the design points $X_{\beta i}$, it may be shown using moment methods that for $\ell = 3$,

(A.13) $$\sum_{x \in \mathcal{R}_k} |S_\ell(x)| = o_p(n_\beta^{-1/2}).$$

Property (3.7) implies that the bias of $\hat{f}_{\beta i}$ is of order $H^{d_1} = O(n_\beta^{-(1/2)-\eta})$ for some $\eta > 0$, whence it may be proved that (A.13) holds with $\ell = 4$. Result (A.11) and the property $n_\beta^{c-1} H^{-k} + H^{2d_1} = O(n_\beta^{-(1/2)-\eta})$ for some $c, \eta > 0$, which follows from (3.7), together imply (A.13) with $\ell = 5$. The lattice argument is used in the cases $\ell = 3, 4, 5$.

Next, we develop approximations to $S_2(x)$. Note that defining $a(v, x) = H^{-k} L\{(x-v)/H\}$, we have

$$\hat{f}_{\beta,-i}(x) = \frac{1}{n_\beta - 1} \sum_{j: j \neq i} a(X_{\beta j}, x).$$

Given $1 \leq i, j \leq n_\beta$ with $i \neq j$, define
$$A(u,v,x) = \frac{g(u)\chi_\beta(u,x)\{a(v,u) - \kappa_\beta(u)\}}{f_\beta(u)^2 h^{k-p}} K\left(\frac{u^{[k-p]} - x^{[k-p]}}{h}\right).$$

We shall construct a $U$-statistic-type projection of $A(X_{\beta i}, X_{\beta j}, x)$ using $D_1(v,x) = E\{A(X_{\beta i}, v, x)\}$, $D_2(u,x) = E\{A(u, X_{\beta j}, x)\}$ and $D_3(x) = E\{A(X_{\beta j}, X_{\beta i}, x)\}$. However, $D_2 \equiv 0$ and therefore $D_3 \equiv 0$, whence $S_2 = T_1 + T_2$, where

$$T_1(x) = \frac{1}{n_\beta} \sum_{i=1}^{n_\beta} D_1(X_{\beta i}, x),$$

$$T_2(x) = \frac{1}{n_\beta(n_\beta - 1)} \sum_{i=1}^{n_\beta} \sum_{j:j \neq i} \{A(X_{\beta i}, X_{\beta j}, x) - D_1(X_{\beta j}, x)\}.$$

Now, $D_1(v,x) = D_3(v,x) - E\{D_3(X_{\beta j}, x)\}$, where

$$D_3(v,x) = E\left[\frac{g(X_{\beta i})\chi_\beta(X_{\beta i}, x) a(v, X_{\beta i})}{f_\beta(X_{\beta i})^2 h^{k-p}} K\left(\frac{X_{\beta i}^{[k-p]} - x^{[k-p]}}{h}\right)\right]$$

$$= E\left[\frac{g(X_{\beta i})\chi_\beta(X_{\beta i}, x)}{f_\beta(X_{\beta i})^2 h^{k-p}} K\left(\frac{X_{\beta i}^{[k-p]} - x^{[k-p]}}{h}\right) L\left(\frac{X_{\beta i} - v}{H}\right)\right].$$

Let $\xi(v,x) = g(v)\chi_\beta(v,x) f_\beta(v)^{-1}$. Then with the $O(n_\beta^{-\eta})$ remainders below being of that form uniformly in $v, x \in \mathcal{R}_k$, for some $\eta > 0$, we have

$$D_3(v,x) = \frac{\xi(v,x) f(v)^{-1} + O(n_\beta^{-\eta})}{h^{k-p} H^k} E\left[K\left(\frac{X_{\beta i}^{[k-p]} - x^{[k-p]}}{h}\right) L\left(\frac{X_{\beta i} - v}{H}\right)\right]$$

$$= \frac{\xi(v,x) + O(n_\beta^{-\eta})}{h^{k-p}} \int K\left(\frac{v^{[k-p]} - x^{[k-p]}}{h} + H h^{-1} w^{[k-p]}\right) L(w)\, dw.$$

(A.14)

Noting that by (3.7), $H h^{-1} = O(n_\beta^{-\eta})$ for some $\eta > 0$ and using the lattice argument, it can be proved from (A.14) that, uniformly in $x \in \mathcal{R}_k$,

(A.15)
$$T_1(x) = \frac{1}{n_\beta} \sum_{i=1}^{n_\beta} (1-E) \frac{\xi(X_{\beta i}, x)}{h^{k-p}} K\left(\frac{X_{\beta i}^{[k-p]} - x^{[k-p]}}{h}\right)$$
$$+ o_p\{(n_\beta h^{k-p})^{-1/2}\},$$

where $E$ denotes the expectation operator. More simply, moment methods and the lattice argument can together be used to show that $T_2(x) = o_p\{(n_\beta h^{k-p})^{-1/2}\}$, uniformly in $x$. This result and (A.15) together imply that, uniformly in $x \in \mathcal{R}_k$,

$$S_2(x) = \frac{1}{n_\beta} \sum_{i=1}^{n_\beta} (1-E) \frac{\xi(X_{\beta i}, x)}{h^{k-p}} K\left(\frac{X_{\beta i}^{[k-p]} - x^{[k-p]}}{h}\right)$$

(A.16)
$$+ o_p\{(n_\beta h^{k-p})^{-1/2}\}.$$

Combining (A.12), (A.13) for $\ell = 3, 4, 5$ and (A.16), we find that, uniformly in $x \in \mathcal{R}_k$,

$$\begin{aligned}
(\widehat{\psi_\beta g^\beta})(x) &= S_1(x) - S_2(x) + o_p\{(n_\beta h^{k-p})^{-1/2}\} \\
&= \frac{1}{n_\beta} \sum_{i=1}^{n_\beta} \frac{\varepsilon_{\beta i} \chi_\beta(X_{\beta i}, x)}{f_\beta(X_{\beta i}) h^{k-p}} K\left(\frac{X_{\beta i}^{[k-p]} - x^{[k-p]}}{h}\right) + E\{S_1(x)\} \\
&\quad + o_p\{(n_\beta h^{k-p})^{-1/2}\}.
\end{aligned} \tag{A.17}$$

[Note that $S_2(x)$ cancels, up to terms of order $o_p\{(n_\beta h^{k-p})^{-1/2}\}$, with $S_1(x) - E\{S_1(x)\}$, except for the part of the latter that involves the errors $\varepsilon_{\beta i}$.] Result (3.10) follows directly from (A.17).

Department of Mathematics and Statistics
The University of Melbourne
Melbourne, Victoria 3010
Australia
E-mail: halpstat@ms.unimelb.edu.au

Department of Economics
University of Toronto
150 St. George Street
Toronto, Ontario
Canada M5S 3G3
E-mail: yatchew@chass.utoronto.ca